\documentclass[a4paper]{jpconf}
\usepackage{graphicx}

\usepackage{amsmath,amsthm}
\usepackage{amssymb}
\usepackage{amsfonts}
\usepackage{amscd}
\usepackage{color}

\newtheorem{thm}{Theorem}[section]

\newtheorem{proposition}[thm]{Proposition}
\theoremstyle{definition}
\newtheorem{definition}[thm]{Definition}

\newcommand{\cR} {\mathcal R}

\newcommand{\id}{\relax{\rm 1\kern-.28em 1}}
\newcommand{\R}{\mathbb{R}}
\newcommand{\C}{\mathbb{C}}

\newcommand{\fp}{\mathfrak{p}}

\newcommand{\re}{\mathrm{e}}
\newcommand{\rspan}{\mathrm{span}}
\newcommand{\rGL}{\mathrm{GL}}

\newcommand{\rSU}{\mathrm{SU}}
\newcommand{\rSL}{\mathrm{SL}}
\newcommand{\rSO}{\mathrm{SO}}

\newcommand{\rM}{\mathrm{M}}

\newcommand{\rP}{\mathrm{P}}

\newcommand{\Gr}{\mathrm{Gr}}
\newcommand{\Fl}{\mathrm{Fl}}

\newcommand{\rGl}{\mathrm{Gl}}

\newcommand{\fsl}{\mathfrak{sl}}
\newcommand{\fsu}{\mathfrak{su}}
\newcommand{\bfour}{\mathbf{4}}

\newcommand{\ri}{{\rm{i}}}
\newcommand{\beq}{\begin{equation}}
\newcommand{\eeq}{{\end{equation}}}

\begin{document}
\title{\centerline{Quantum Chiral Superfields}}

\author{\centerline{
Rita Fioresi$^*$,
Mar\'{\i}a A. Lled\'{o}$^\dagger$,
Junaid Razzaq$^\bullet$}}

\address{\centerline{$^*$
    FaBiT, via San Donato 15, Universit\`{a} di Bologna, Bologna, Italy}}
\address{\centerline{$^\dagger$Department de F\'{\i}sica Te\`{o}rica, Universitat de Val\`{e}ncia and IFIC (CSIC).}
 \centerline{Carrer Dr. Moliner, 50. 46100 Burjassot. Spain.}
}
\address{\centerline{$^\bullet$ Dipartimento di Matematica,
    Piazza Porta San Donato 5,
    Universit\`{a} di Bologna, Italy}}

\ead{$^*$rita.fioresi@unibo.it, $^\dagger$maria.lledo@ific.uv.es, $^\bullet$junaid.razzaq2@unibo.it}

\begin{abstract}
  We define the ordinary Minkowski space inside the conformal
  space according to Penrose and Manin as homogeneous spaces
  for the Poincar\'e and conformal group respectively.
  We realize the supersymmetric
  (SUSY) generalizations of such homogeneous spaces over the complex
  and the real fields.
  We finally investigate chiral (antichiral) superfields,
  which are superfields on
  the super Grassmannian, $\Gr(2|1, 4|1)$, respectively
  on $\Gr(2|0,4|1)$. They ultimately give the twistor coordinates necessary
  to describe the conformal superspace as the
  flag $\mathrm{Fl}(2|0,2|1;4|1)$ and the Minkowski superspace
  as its big cell.
\end{abstract}

\section{Introduction}

Quantum Groups are born to encode quantum symmetries. In particular, to be
able to treat the quantization of spacetime symmetries. Physically significant geometric objects as
the {\sl conformal} or the  {\sl Minkowski space} are
 homogeneous spaces, which have a natural action of the conformal and Poincar\'{e} groups respectively. In the language of non commutative
geometry, it  is natural to encode the quantization of any geometric
theory, as for example, gravity or Yang-Mills theories, in terms  of algebraic and differential geometry.

Our spacetime can be, ultimately, a {\sl superspace}  with even and odd coordinates. Field theories having {\sl supersymmetry} have symmetries that can mix even and odd coordinates and predict that each particle must have a {\it superpartner}, with the same mass but with opposite character (one is bosonic, the other is fermionic). At least for $N=1$ supersymmetry, supersymmetric theories can be formulated over $N=1$ superspace.

We have successfully
quantized, as homogeneous spaces, the {\sl Minkowski}
and {\sl conformal superspaces}, together with the  coactions of their natural symmetries:
the  {\sl  conformal} and {\sl Poincar\'{e} quantum supergroups}.

However, a new emerging theory of
{\sl Quantum Cartan Geometry} holds the potential
to successfully go towards a quantum theory of
  Cartan connections and (bi)covariant objects (e.g.
  covariant hamiltonian) and quantum (super) differential calculus. Then one could, in principle, attack  the far reaching application to
{\sl non commutative (super) gravity}.

\medskip
Our paper is organized as follows.

In Section \ref{ordinary} we recall the construction of the ordinary Minkowski
space inside the conformal space as proposed by Penrose \cite{pe} and Manin \cite{ma}: the complex, conformal space
is the {\sl Grassmannian} of two spaces in a four dimensional space $\Gr(2,4)$ and Minkowski space is its {\sl big cell}.

In Section \ref{minkowskisuper} we give the complex conformal superspace as the {\sl flag supermanifold } $\Fl(2|0,2|1,4|1)$, and the   complex Minkowski superspace as its big cell \cite{ma}.  We then compute the relevant  real forms of these objects. As in the non super case, the structure of superflag disappears when imposing the reality condition. The superfields in this real Minkowski superspace are  relevant in physics for the description of the electromagnetic or Yang-Mills theories, since they contain, as the highest spin particle, a vector field (in fact, they are called vector superfields). We also investigate chiral (antichiral) superfields, which are superfields on the superspace of a super Grassmannian, $\Gr(2|1,4|1)$, ($\Gr(2|0,4|1)$). These superspaces are genuinely complex, since a real condition is not be compatible with the action of the Lorentz group.  They contain a spin 0 and a spin 1/2 particle. The theory of the interacting chiral and vector superfields is the supersymmetric version of QED. Also, the supersymmetric gauge transformations are described in terms of chiral superfields.

In Section \ref{quantumsuper}, we proceed to give a quantum deformation of the
geometric objects that we have studied so far in $N=1$ supersymmetry.

This paper is a summary of results obtained in a series of papers by the authors, see for example \cite{cfl,fl,flln,flr}.

\section{The ordinary Minkoski space}\label{ordinary}

We start by giving the realization of the conformal algebra in Minkowski space time. Let $\{x^\mu\}_{\mu=0}^4$ be the coordinates of Minkowski space.  Then we have the following generators of the conformal transformations

\begin{align*}&P_\mu=\partial_\mu\,,\qquad &\hbox{(translations)},\\
&D= x^\mu\partial_\mu\,,\qquad &\hbox{(dilation)},\\
&L_{\mu\nu}= x_\nu\partial_\mu-x_\mu\partial_\nu\,,\qquad &\hbox{(Lorentz transformations)},\\
&K_\mu=2x_\mu x^\nu\partial_\nu-x^2\partial_\mu\,,\qquad &\hbox{(special conformal transformations)}.
\end{align*}
The conformal (finite) transformation on Minkowski space are:
\begin{align*}&  {x'}^\mu= x^\mu+ a^\mu\,,\qquad &\hbox{(translations)},\\
& {x'}^\mu=\re^{u/ 2}x^\mu\,,\qquad &\hbox{(dilation)},\\
& {x'}^\mu=\Lambda^\mu_{\phantom{\mu}\nu} x^\nu\,,\qquad &\hbox{(Lorentz transformations)},\\
& {x'}^\mu= \frac{x^\mu+b^\mu x^2}{1+ 2 b\cdot x+2b^2x^2} \,,
\qquad &\hbox{(special conformal transformations)}\,.
\end{align*}
The quantities $a^\mu, u, b_\mu$ plus the six parameters of the Lorentz group in $\Lambda^\mu_{\phantom{\mu}\nu}$ are the parameters of the conformal transformations. We notice that, for fixed values of $b_\mu$ the special conformal transformations are singular at some points of spacetime. This suggests that the Minkowski space is not the complete conformal space, were the conformal group is realized, and that some points have to be added to deal with these infinities.

The conformal group is in fact $\rSO(2,4)$, and its spin group (the double cover) is $\rSU(2,2)$. We will work with the last version, which is needed when we introduce fermions, as in supersymmetry.

\subsection{The ordinary complex Minkowski space as the big cell
in the Grassmannian $G(2,4)$}

The approach of Penrose and Manin is done in the complexification of spacetime. The complexification of the conformal group (above) is $\rSL(4,\C)$. We consider $\C^4$ (not Minkowski spacetime) and the Grassmannian of two planes inside $\C^4$, $\Gr(2,4)$. The Grassmannian is an homogeneous space for $\rSL(4,\C)$, the isotropy group being the parabolic subgroup
\begin{equation} Q=\left\{\,\begin{pmatrix}* & * \\0 & *\end{pmatrix}\,\right\}, \qquad \Gr(2,4) = \rSL_4(\C)/Q\,.\label{isotropy}\end{equation}  An element of the Grassmannian is given by two linearly independent vectors
$$\rspan\left\{\begin{pmatrix}a_1\\b_1\\c_1\\d_1\end{pmatrix}, \begin{pmatrix}a_2\\b_2\\c_2\\d_2\end{pmatrix}\right\}\,.$$
Then, some of the minors of the matrix formed by the two column vectors must be non zero. Let us suppose that the minor
\begin{equation}\det\begin{pmatrix}
a_1 & a_2 \\
   b_1 & b_2
\end{pmatrix}\neq 0\label{determinant}\,,\end{equation}
then, there is a
linear combination of these vectors such that the $4\times 2$ matrix is of the form
$$\begin{pmatrix}
    1 & 0 \\
   0 & 1\\
   t_1& t_2\\
   t_3&t_4
 \end{pmatrix}\,.$$
 We are using here the right action of $\rGL(2, \C)$ over this set. The open set where the determinant (\ref{determinant}) is different from zero is the big cell of the Grassmannian. It is the complex affine space of dimension 4. Its relation with the Minkowski space is through  the Pauli matrices

 $$\sigma_0=\begin{pmatrix}
   1& 0 \\
  0 & 1
 \end{pmatrix},\qquad \sigma_1=\begin{pmatrix}
   0& 1 \\
  1 & 0
 \end{pmatrix},\qquad\sigma_2=\begin{pmatrix}
   0& -\ri \\
  \ri & 0
 \end{pmatrix},\qquad\sigma_3=\begin{pmatrix}
   1& 0 \\
  0 & -1
 \end{pmatrix}\,,$$
 with

$$A=\begin{pmatrix}
    t_1& t_2\\
   t_3&t_4
 \end{pmatrix}=x^\mu\sigma_\mu=\begin{pmatrix}x^0+x^3&x^1-\ri x^2\\x^1+\ri x^2& x^0-x^3
\end{pmatrix},\qquad \mu=0, 1,2,3,4$$
 and sum over repeated indices is understood. The coordinates $x^\mu$ are the standard coordinates of the Minnkowski space. This is justified if we see how the group $\rSL(4,\C)$ acts on $A$. In fact, the subgroup that leaves invariant the bigcell, and its action on it is

$$
\begin{pmatrix}L& 0 \\NL & R\end{pmatrix}
\begin{pmatrix}I \\ A \end{pmatrix}Q\,= \,
\begin{pmatrix}L \\ NL+RA \end{pmatrix}Q \, = \,
\begin{pmatrix}I \\ N+RAL^{-1} \end{pmatrix}Q \,.
$$

where
\begin{equation}
 N=n^\mu\sigma_\mu=
\begin{pmatrix}n^0+n^3&n^1-\ri n^2\\n^1+\ri n^2& n^0-n^3\end{pmatrix}\,.
\label{ordinarycoordinates}\end{equation}
This subgroup is the complexified Poincar\'{e} group plus dilations, $N$ representing the translations, $L$ and $R$ in $\rGL(2,\C)$ with $\det L\cdot \det R=1$.

%\subsection{Penrose approach: the real Minkowski space}
%The real form of $\rSL(4,\C)$  gives for the Poincar\'{e} group the subgroup   $P=\rSL_2(\C) \otimes \R^4$ (remember that $\fsl_2(\C) \cong \fso(3,1)$).
%
%

The Minkowski metric
$$\det(A)=(x^0)^2-(x^1)^2-(x^2)^2-(x^3)^2$$
is automatically preserved by Poincar\'{e} action.

\bigskip

The  real forms of these objects, consistent with the Minkowski signature, are well known. The real  conformal group is $\rSU(2,2)$ (the spin group or  double covering of $\rSO(2,4)$). The conformal space is not anymore a Grassmannian, but the conformal compactification of the Minkowski space. The real Miknowski space has all the coordinates $x^\mu$ in (\ref{ordinarycoordinates}) real, so $A$ is hermitian. The real  Poincar\'{e} group has $L={R^\dagger}^{-1}$ and $N$ hermitian (as $A$, since they are the translations).

\section{Minkowski superspace}\label{minkowskisuper}

\medskip
\subsection{Wess-Zumino superalgebra: $N=1$ supersymmetry}

The Wess-Zumino complex conformal superalgebra ($N=1$) is
$$
\fsl(4|1)=\Biggl\{\begin{pmatrix}p&\alpha\\\beta&c\end{pmatrix}
\;\;|\;\;p\in \rM_4(\C),\; \alpha,\beta^t\in \C^4,\; c\in \C,\; \tr \,
p=c\Biggr\}
$$
with
$$
\fsl(4|1)_0=\Biggl\{\begin{pmatrix}p&0\\0&c\end{pmatrix}\Biggr\}=
\fsl_4(\C)\oplus \C\,,\quad \fsl(4|1)_1=
\Biggl\{\begin{pmatrix}0&\alpha\\\beta&0\end{pmatrix}\Biggr\}=
\bfour^{-1}\oplus {\mathbf{\overline {4}}^{+1}}\,.$$
We observe that

\begin{itemize}
\item
The even part $\fsl(4|1)$ contains, not only the conformal algebra,
but also an \textit{inevitable} extra factor $\C$.
\item The odd part is a spinorial representation
\end{itemize}

The superconformal group is\footnote{Using the functor of points formalism, the entries of this matrix are even or odd elements of an arbitrary superring. We will understand this for all the supergroups and homogeneous spaces without further mention, unless needed.}
$$\rSL(4|1)(\cR)=\left\{ g=\begin{pmatrix}D&\tau\\\gamma&d\end{pmatrix}\; |\; D\in \rGL(4)(\cR),\;\; \tau,\gamma^t\in \cR^4,\;\; d\in \cR\, \mathrm{invertible},\;\; \mathrm{Ber}g=1\right\}$$ where $\cR$ is an arbitrary superring and $\mathrm{Ber}g$ is the Berezinian of $g$.

\medskip
The real form of the complex conformal superalgebra, $\fsl(4|1)$,  that, as we will see, is compatible with the signature of Minkowski space is

\begin{align}\fsu(2,2|1)=\left\{\begin{pmatrix}p&\alpha\\\beta&\ri z\quad
\end{pmatrix}
\in \rM_4(\C)\; | \quad \alpha, \beta^t\in \C^4,\quad z\in \R;
\right.\nonumber\\
\left. Fp+p^\dagger F=0,\;    \tr \, p =\ri z,\; \alpha=\ri F\beta^\dagger\right\}\,,\label{su22}\end{align}

where

$$F=\ri \begin{pmatrix}0&-\id_2\\\id_2&0\end{pmatrix}\,.$$
It can be seen as the space of fixed points of  the conjugation $\sigma$.
 $$\begin{CD}
\fsl(4|1)@>\sigma>>\fsl(4|1)\\
\begin{pmatrix} p&\alpha\\ \beta&d\\ \end{pmatrix} @>>>\begin{pmatrix}
-Fp^\dagger F& \ri F\beta^\dagger\\ \ri \alpha^\dagger F&- {\bar d}\\
\end{pmatrix}
\end{CD}
$$

So $\fsu(2,2|1)=\fsl(4|1)^\sigma$.

The conjugation $\rho$ whose set of fixed points defines the real form $SU(2,2|1)$ is

\begin{equation}\begin{CD}\rSL(4|1)@>\rho>>\rSL(4|1)\\
\begin{pmatrix}D&\tau\\
\gamma&d\end{pmatrix}@>>>L\begin{pmatrix}D^\dagger&i\gamma^\dagger\\\ri\tau^\dagger&\bar d\end{pmatrix}^{-1}L\,
\end{CD}$$
with
$$L=\begin{pmatrix}F&0\\0&1\end{pmatrix}\,.\label{realform}\end{equation}

It can be proven that this conjugation induces in the Lie algebra $\fsl(4|1)$ the conjugation (\ref{su22}).

All these real forms are explicitly computed in \cite{fl}.

\subsection{Minkowski superspace and the super Poincar\'{e} group}
The $N=1$ Minkowski superspace is an affine superspace that, when complexified, can be seen as the big cell of a flag supermanifold. The flag supermanifold is then the conformal superspace.
We first describe it and then we will deduce the properties of the Minkowski superspace as the big cell inside it.

So we consider the flag supermanifold $\Fl(2|0,2|1, 4|1)$. It is an homogenous space of the supergroup $\rSL(4|1)$, with isotropy group
$$F_u=\left\{\begin{pmatrix}
g_{11}&g_{12}&g_{13}&g_{14}&\gamma_{15}\\
g_{21}&g_{22}&g_{23}&g_{24}&\gamma_{25}\\
0&0&g_{33}&g_{34}&0\\
0&0&g_{43}&g_{44}&0\\
0&0&\gamma_{53}&\gamma_{54}&g_{55}\\
\end{pmatrix}\right\}\subset\rSL(4|1)\,.$$ This is similar to (\ref{isotropy}) in the non super case.  An element of the superflag is given by two sets of vectors, one describing a $2|0$ subspace, so it is the span of two linearly independent even vectors.
$$P_1=\begin{pmatrix}
a_{11}&a_{12}\\
a_{21}&a_{22}\\
a_{31}&a_{32}\\
a_{41}&a_{22}\\
\gamma_{51}&\gamma_{52}
\end{pmatrix}$$
and the other  is the span of two even and one odd vector

$$P_2=\begin{pmatrix}
b_{11}&b_{12}&\delta_{13}\\
b_{21}&b_{22}&\delta_{23}\\
b_{31}&b_{32}&\delta_{33}\\
b_{41}&b_{42}&\delta_{43}\\
\delta_{51}&\delta_{52}&b_{53}\\
\end{pmatrix}\,.$$
Then we will have to impose the condition that $P_1\subset P_2$. This condition will result in the {\sl twistor} or {\sl incidence relations}.

As for the non super case, we assume that
$$\det\begin{pmatrix}
a_{11}&a_{12}\\
a_{21}&a_{22}\end{pmatrix},\quad \det\begin{pmatrix}
b_{11}&b_{12}\\
b_{21}&b_{22}\end{pmatrix}\quad\hbox{and}\quad b_{53} $$
are invertible
The open set were these conditions are satisfied is the big cell of the flag supermanifold. For $P_1$ we can use the right action of $\rGl(2)$ to put it in standard form

$$
\begin{pmatrix}
\id \\
A \\
\alpha
\end{pmatrix}\,.$$
For $P_2$ we use the right action of $\rGL(2|1)$ and the result is

$$\begin{pmatrix}
\id & 0 \\
B & \beta \\
0 & 1
\end{pmatrix}\,.$$
The twistor relations are
$$B=A-\beta\alpha\,,$$
so $(A,\alpha,\beta)$ or $(B,\alpha, \beta)$ are both sets of global coordinates on the complex Minkowski superspace. The super Poincar\'{e} group is the super group that leaves the big cell invariant

$$\rP(\cR)=\left\{\begin{pmatrix}L&0&0\\NL&R&R\varphi\\\chi&0&d\end{pmatrix}\;|\;\; L, R\in \rGL_2(\cR), \;N\in \rM_2(\cR),\; \chi, \varphi\in\cR^2,\;d\in\cR \, \mathrm{invertible}\right\}$$

The complex Poincar\'{e} superalgebra is

$$\mathfrak{p}_\C=\left\{\begin{pmatrix}l&0&0\\n&r&\rho\\\delta&0&f\end{pmatrix}\;|\;\; l,r,n\in \rM_2(\C), \; \rho, \delta\in\C^2,\;f\in\C\right\}$$

The real Poincar\'{e} superalgebra is
$$
\fp_{\R}=\left\{\begin{pmatrix} l&0&0\\ n&-l^\dagger&-\delta^\dagger \\ \delta&0&\ri z
\end{pmatrix}\,|\; n, l\in M_2(\C), \delta\in \C^2, z\in \R \right\}\,.$$

 For the Poincar\'{e} supergroup the conjugation can be computed more explicitly:
\begin{equation}\begin{CD}\begin{pmatrix}L&0&0\\M&R&R\varphi\\d\chi&0&d\end{pmatrix}@>>>\begin{pmatrix}
{R^\dagger}^{-1}&0&0\\
{L^\dagger}^{-1}M^\dagger{R^\dagger}^{-1}+{L^\dagger}^{-1}\chi^\dagger\varphi^\dagger&{L^\dagger}^{-1}&{L^\dagger}^{-1}\chi^\dagger\\
{\bar d}^{-1}\varphi^\dagger&0&{\bar d}^{-1}
\end{pmatrix}\,,
\end{CD}\label{poincaresupergroup}\end{equation}
so the reality condition is
$$L={R^\dagger}^{-1},\qquad \varphi=\chi^\dagger,\qquad ML^{-1}={ML^{-1}}^\dagger+{L^\dagger}^{-1}\chi^\dagger\chi L^{-1}\,.$$
With the definitions
$$N:=ML^{-1},\qquad T:=N-\frac 12{L^\dagger}^{-1}\chi^\dagger\chi L^{-1}\,,$$
the reality condition becomes
$$L={R^\dagger}^{-1},\qquad \varphi=\chi^\dagger,\qquad T=T^\dagger\,.$$

\medskip

The flag supermanifold is an homogeneous space  ($\rSL(4|1/F_u$). Taking the coordinates $(A, \alpha,\beta)$, a coset representative is
$$g=\begin{pmatrix}\id&0&0\\A&\id&\beta\\\alpha&0&1\end{pmatrix}\,.$$ From the real form (\ref{realform}) one obtains
$$A=A^\dagger+\alpha\alpha^\dagger,\qquad \beta=\alpha^\dagger\,.$$
Defining
$$C:=A-\frac 12\beta\alpha$$
the reality condition becomes
$$C=C^\dagger,\qquad \beta=\alpha^\dagger\,.$$

We change now the notation to the one used in physics:
$$(\beta,\alpha)\longrightarrow (\theta,{\bar\theta}^t)\,,$$ then the big cell can be expressed as

$$\left\{\begin{pmatrix}
C+\frac 12 \theta\bar\theta^t \\
\bar \theta^t
\end{pmatrix}, \,
\begin{pmatrix}
\id & 0 \\
C -\frac 12 \theta\bar\theta^t  & \theta \\
0 & 1
\end{pmatrix} \right\}\,.$$
Then the action of the super Poincar\'{e} group on the big cell is

\begin{align}
 &C\;\longrightarrow \;R\bigl(C+\frac 12\varphi\bar\theta^t-\frac 12 \theta\bar\varphi^t\bigr)R^\dagger+T,\nonumber\\
 &\theta\;\longrightarrow\; d^{-1}R(\theta+\varphi),\nonumber\\
 &\bar \theta\;\longrightarrow\; dL^{-1 t}(\bar\theta+\bar \varphi)
\label{actionchanged}\,.\end{align}

We retrieve the usual Minkowski coordinates as
$$
C=\sum_{\mu=0}^3x^\mu\sigma_\mu=\begin{pmatrix}x^0+x^3&x^1-
\ri x^2\\x^1+\ri x^2& x^0-x^3\end{pmatrix}\,.
$$
Notice that the superflag is embedded in
$$\Fl(2|0,2|1,4|1)\subset \Gr(2|0,4|1)\times \Gr(2|1,4|1)\,.$$
 The big cell of $\Gr(2|1,4|1)$ is the chiral Minkowski superspace. Its superfields
 depend only on $x^\mu$ and $\theta$. In the same way, the antichiral Minkowski superspace is the big cell of $\Gr(2|0,4|1)$, and superfields depend only on $x^\mu$ and $\bar\theta$.i

This treatment can be extended to $N=2$ supersymmetry \cite{flr}.

%\subsection{Chiral Superfields $N=1$ ($N=2$)}
%
%{\it chiral superfield} is a superfield $\Phi$ such that
%$$
%\bar D_{\dot\alpha}\Phi=0
%$$
%where  $D_{\dot\alpha}$ are the infinitesimal generators
%of the supertranslation superalgebra in Poincar\`e:
%$$
% P_\mu=\partial_\mu,\qquad D_\alpha=\nabla_\alpha-\frac 14\bar\theta^{\dot\alpha}\bar\sigma^\mu_{\dot\alpha\alpha}\partial_\mu, \qquad
%\bar D_{\dot \alpha}=\bar\nabla_{\dot\alpha}-\frac 14\bar\sigma^\mu_{\dot\alpha\alpha}\theta^{\alpha}\partial_\mu
%$$
%with commutation relations:
%$$
%[D_\alpha,\bar D_{\dot \alpha}]
%=-\frac 12 \bar\sigma^\mu_{\dot\alpha\alpha}P_\mu,\qquad \hbox{and the rest zero}
%$$
%The supertranslation
%algebra acts on the space of chiral superfields (chiral subspace)
%and the whole Poincar\`e does.
%
%
%
%

\section{Quantum Minkowski Superspace}\label{quantumsuper}
% as Homogeneous
%space, $N=1$, ($N=2$)}

The  key idea is to replace the geometric objects with their
non commutative function algebras of {\sl polynomials}. Let
$\C_q:=\C[q,q^{-1}]$ be Laurent polynomials in $q=e^h$. We have for example, the replacements:

\begin{align*}
&M=\C^{4}
\quad \longrightarrow \quad \C_q[M],  \\
&P=\left(\rSL_2(\C) \times \rSL_2(\C) \times\C\right)\ltimes \C^4
\longrightarrow \quad \C_q[P]\,, %\quad \hbox{ring with gen. and rels.}
\end{align*}
where $\C_q[M]$ and $\C_q[P]$ are rings given in terms of generators and relations.
One also replaces actions with {\sl coactions}:
$$
P\times M\rightarrow M\quad \longrightarrow \quad\C_q[M]\rightarrow\C_q[P] \otimes \C_q[M]\,.
$$
The same can be done in the super case, where we  will end up with non commutative superalgebras. In the case of the Lie groups and supergroups appearing here, the polynomial algebras have the structure of Hopf (super)algebras

\bigskip

Up to now we have studied the real Minkowski superspace. For the quantization, we take a simpler problem, that is, the chiral (or antichiral: they are both isomorphic) superspace. The difference is significant: For the quantum superflag the commutation relations among the generators are much more difficult to compute, since they describe the two dual Grassmannians $\Gr(2|0,4|1)$ and $\Gr(2|1,4|1)$ and, moreover, one would have to compute the quantum twistor relations. In fact, the problem of quantizing the flag supermanifold  $\Fl(2|0,2|1,4|1)$ was approached in \cite{flln} with a less direct formalism, in terms of a quantum section.

So our aim is to quantize $\Gr(2|0,4|1)$ and its big cell substituting $\rSL(4|1)$ by $\rSL_q(4|1)$. For the quantum supergroup, we  take Manin's approach \cite{ma2}. In a unified notation, an element of $\rSL_q(4|1)$ is denoted as $[a_{ij}]$ with $i,j=1,\dots, 5$. The entries $a_{i5}, a_{5i}$ with $i=1,\dots , 4$ are odd elements of the arbitrary superring $\cR$ and  the rest, even elements. Then, the commutation relations  are

\begin{align*}
&a_{ij}a_{il}=(-1)^{\pi(a_{ij})\pi(a_{il})}
q^{(-1)^{p(i)+1}}a_{il}a_{ij}, && \hbox{for  } j < l \\&&& \\
&a_{ij}a_{kj}=(-1)^{\pi(a_{ij})\pi(a_{kj})}
q^{(-1)^{p(j)+1}}a_{kj}a_{ij}, && \hbox{for  } i < k \\ &&&\\
&a_{ij}a_{kl}=(-1)^{\pi(a_{ij})\pi(a_{kl})}a_{kl}a_{ij}, &&  \hbox{for  }
i< k,\;j > l \\&&&\hbox{or } \quad i > k,\; j < l \\&&& \\
&a_{ij}a_{kl}-(-1)^{\pi(a_{ij})\pi(a_{kl})}a_{kl}a_{ij}=(-1)^{\pi(a_{ij})\pi(a_{kl})}(q^{-1}-q)
a_{kj}a_{il},&&\\
&&& \hbox{for  }\quad i<k,\;j<l\,,
\end{align*}
where we have assigned a parity to the indices: $p(i)=0$ if $1 \leq i \leq 4$ and  $p(5)=1$. The parity of $a_{ij}$ is  $\pi(a_{ij})=p(i)+p(j)$ mod 2. It is shown in  \cite{cfl} (see also \cite{fl} for a more detailed exposition) that one can define, mimicking what happens in the non quantum case:

\begin{definition}The quantum super Grassmannian $\Gr_q:=Gr_q(2|0,4|1)$ is  the subalgebra of $\rSL_q(4|1)$ generated by the quantum minors\label{qgrassmannian-def}
$$
D_{ij}:=a_{i1}a_{j2}-q^{-1}a_{i2}a_{j1}\qquad D_{i5}:=a_{i1}a_{52}-q^{-1}a_{i2}a_{51}\qquad
D_{55}:=a_{51}a_{52}
$$
with $1\leq i<j\leq 4$.
\end{definition}
The point here is to show that the commutation relations close over these elements. In fact, one can prove that this is the  case. They are given explicitly in \cite{cfl,fl}. The quantum supergroup has a coaction over $\Gr_q$

$$\begin{CD}
 \Gr_q@> \Delta>> \rSL_q(4|1)\otimes\Gr_q\,.
 \end{CD}
 $$
that follows by restricting the comultiplication of the quantum supergroup.

What about the quantum Minkowski superspace? We have to restrict to the big cell, which means, in the quantum setting, to localize the algebra by inverting one of the quantum minors, in particular, $D_{12}$. This localization is denoted as $\Gr_q[D_{12}^{-1}]$.

In the classical case, the group that leaves invariant the quantum big cell of the chiral Minkowski space is bigger than the super Poincar\'{e} group. It is a lower parabolic
subgroup of  $\rSL(4|1)$.
 \begin{equation}
\begin{pmatrix}x&0&0\\
 tx&y& y\eta\\ \tau x&\xi&d\end{pmatrix}\,,\label{qps} \end{equation}

Taking module by $\xi$ we obtain the super Poincar\'{e} group. In the quantum case, formally, the quantum super Poincar\'{e} group looks like (\ref{qps}) (modulo $\xi$). But now the entries satisfy the commutation relations of the Manin supermatrix:

  \begin{align*}
&x =\begin{pmatrix}
a_{11} & a_{12}\\ a_{21} & a_{22}  \end{pmatrix},
&&
t =
\begin{pmatrix} -q^{-1}D_{23}D_{12}^{-1} & D_{13}D_{12}^{-1}\\
-q^{-1}D_{24}D_{12}^{-1} & D_{14}D_{12}^{-1}  \end{pmatrix},
\nonumber \\\nonumber \\
&y = \begin{pmatrix}
a_{33} & a_{34}\\ a_{43} & a_{44} \end{pmatrix},&&
d =a_{55}, \\ \\
&
\tau=(  -q^{-1}D_{25}D_{12}^{-1},
D_{15}D_{12}^{-1} ),&& \eta =
\begin{pmatrix}
-q^{-1}{D^{34}_{34}}^{-1}D^{45}_{34} \\
{D^{34}_{34}}^{-1}D_{34}^{35} \\
\end{pmatrix}\,.
\end{align*}

\begin{definition}
We define the quantum chiral Minkowski superspace as the subalgebra of the quantum super Poincar\'{e} group generated by

$$t =
\begin{pmatrix} -q^{-1}D_{23}D_{12}^{-1} & D_{13}D_{12}^{-1}\\
-q^{-1}D_{24}D_{12}^{-1} & D_{14}D_{12}^{-1}  \end{pmatrix}, \qquad \tau=(  -q^{-1}D_{25}D_{12}^{-1},
D_{15}D_{12}^{-1} )\,.
$$

\hfill$\blacksquare$

\end{definition}

Equivalently, it is a computation to prove

\begin{proposition} \label{commt-tau}
The quantum chiral Minkowski superspace has
the following  presentation:
$$
\C_q[\rM] = \C_q \langle t_{ij}, \tau_{5j} \rangle \, \big/ \, I_U,
\qquad 3 \leq i \leq 4, \; j=1,2\,,
$$
where $I_U$ is the ideal generated by the relations:
\begin{align*}
&t_{i1}t_{i2}=q \, t_{i2}t_{i1}, \qquad
&&t_{3j}t_{4j}=q^{-1} \, t_{4j}t_{3j},
\qquad 1 \leq j \leq 2, \quad 3 \leq i \leq 4
\\ \\
&t_{31}t_{42}=t_{42}t_{31}, \qquad
&&t_{32}t_{41}=t_{41}t_{32}+(q^{-1}-q)t_{42}t_{31},
\\ \\
&\tau_{51}\tau_{52}=-q^{-1}\tau_{52}\tau_{51}, \qquad
&&t_{ij}\tau_{5j}=q^{-1}\tau_{5j}t_{ij}, \qquad 1 \leq j \leq 2 \\ \\
&t_{i1}\tau_{52}= \tau_{52}t_{i1},
\qquad  &&t_{i2}\tau_{51}=\tau_{51} t_{i2}+(q^{-1}-q)t_{i1}\tau_{52}\,.
\end{align*}

\hfill$\blacksquare$
\end{proposition}

This quantum superspace admits a coaction of the quantum super Poincar\'{e} group, computed by restricting the comultiplication in the quantum supergroup.

\begin{align*}&\Delta(t)=t\otimes 1+yS(x)\otimes t+y\eta S(x)\otimes\tau,
&\Delta(\tau)=d\otimes 1)(\tau\otimes 1+1\otimes \tau)(S(x)\otimes 1\,,
\end{align*}
where $S$ denotes the antipode in the group.

\section*{ Acknowledgements}

This work is supported by the Spanish Grant PID2020-116567GB- C21 funded by MCIN/AEI/10.13039/501100011033, by the project PROMETEO/2020/079 (Generalitat Valenciana, by HORIZON-MSCA-2022-SE-01-01 CaLIGOLA and by  Gnsaga-Indam. 

This publication is based upon work from COST Action CaLISTA CA21109 supported by COST (European Cooperation in Science and Technology). www.cost.eu.

\end{document}